\documentclass{amsart}

\usepackage{amsmath,amssymb,amscd,amsfonts}

\newtheorem{theorem}{Theorem}
\newcommand{\bt}{\begin{theorem}}
\newcommand{\et}{\end{theorem}}
\newtheorem{lemma}{Lemma}
\newcommand{\bl}{\begin{lemma}}
\newcommand{\el}{\end{lemma}}
\newtheorem{corollary}{Corollary}
\newcommand{\bc}{\begin{corollary}}
\newcommand{\ec}{\end{corollary}}
\newcommand{\beq}{\begin{equation}}
\newcommand{\eeq}{\end{equation}}
\newcommand{\benum}{\begin{enumerate}}
\newcommand{\eenum}{\end{enumerate}}

\begin{document}

\title{On sequences without geometric progressions}
\author{Melvyn B. Nathanson}
\address{Department of Mathematics\\
Lehman College (CUNY)\\
Bronx, NY 10468}
\email{melvyn.nathanson@lehman.cuny.edu}

\author{Kevin O'Bryant}
\address{Department of Mathematics\\
College of Staten Island (CUNY)\\
Staten Island, NY 10314}
\email{obryant@mail.csi.cuny.edu}


\date{\today}

\begin{abstract}
An improved upper bound  is obtained for the density of sequences
of positive integers that contain no $k$-term geometric progression.
\end{abstract}

\maketitle

\section{A problem of Rankin}

Let $k \geq 3$ be an integer.  Let $r \neq  0, \pm1$ be a real number.
A \emph{geometric progression of length $k$ with common ratio $r$ } is a sequence
$(a_0, a_1,a_2,\ldots, a_{k-1})$ of nonzero real numbers such that
\[
r = \frac{a_i}{a_{i-1}}
\]
for $1, 2, \ldots, k-1$.
For example, $(3/4, 3/2,3,6)$ and $(8, 12, 18,27)$ are geometric progressions
of length 4 with common ratios 2 and  $3/2$, respectively.
A  \emph{$k$-geometric progression} is a geometric progression of length $k$
with common ratio $r$  for some $r$.
If the sequence $(a_0, a_1,a_2,\ldots, a_{k-1})$ is a $k$-geometric progression,
then $a_i \neq a_j$ for $0 \leq i < j \leq k-1$.

A finite or infinite set  of real numbers is \emph{$k$-geometric progression free}
if  the set does not contain numbers $a_0,a_1,\ldots, a_{k-1}$ such that
the sequence $(a_0,a_1,\ldots, a_{k-1})$  is a $k$-geometric progression.
Rankin~\cite{rank60} introduced $k$-geometric progression free sets,
and proved that there exist infinite $k$-geometric progression free sets 
with positive asymptotic density.\footnote{If $A(n)$ 
denotes the number of positive integers 
$a\in A$ with $a \leq n$, then the \emph{upper asymptotic density} of $A$ is 
$d_U(A) = \limsup_{n\rightarrow\infty}A(n)/n$, and 
the  \emph{asymptotic density} of $A$ is 
$d(A) = \lim_{n\rightarrow\infty}A(n)/n$, if this limit exists.}
For example, the set $\mathit{Q}$ of square-free positive integers, 
with asymptotic density $d(\mathit{Q}) = \pi^2/6$,
contains no $k$-term geometric progression for $k \geq 3$.

Let $A$ be a set of positive integers that contains no $k$-term geometric progression.  
Brown and Gordon~\cite{brow-gord96} proved\footnote{Brown and Gordon 
claimed a slightly stronger result, but their proof contains an (easily corrected) error.}
that  the upper asymptotic density of $A$, denoted $d_U(A)$, 
 has the following upper bound:  
\[
d_U(A) \leq 1- \frac{1}{2^k} - \frac{2}{5}\left( \frac{1}{5^{k-1} } -  \frac{1}{6^{k-1}} \right).
\]
Riddell~\cite{ridd69}  and 
Beiglb{\"o}ck, Bergelson, Hindman,  and Strauss\cite{beig-berg-hind-stra06}
proved that 
\[
d_U(A) \leq 1 -  \frac{1}{2^k-1}.  
\]
The purpose of this note is to improve these  results.

\section{An upper bound for  sets with no $k$-term geometric progression}

\bt     \label{GPF:theorem}
For integers $k \geq 3$ and $n \geq 2^{k-1}$, let $GPF_k(n)$ denote the set of subsets
of $\{1,2,\ldots, n\}$ that contain no $k$-term geometric progression.
If $A \in GPF_k(n)$, then 
\[
n - |A| \geq  
\left( \frac{1}{2^k -1} + \frac{2}{5}\left( \frac{1}{5^{k-1} } -  \frac{1}{6^{k-1}} \right)
+  \frac{4}{15}\left( \frac{1}{7^{k-1} } -  \frac{1}{10^{k-1}} \right)  \right)  n
+ O\left(\frac{\log n}{k} \right).
\]
\et

\begin{proof}
Let
\[
L = \left[ \frac {\log 2n}{k \log 2}\right].
\]
For  $1 \leq \ell \leq L$
we have $2^{\ell k -1} \leq n.$
Let $a$ be an odd positive integer such that
\[
a \leq \frac{n}{2^{\ell k -1}}.
\]
The sequence
\[
\left( 2^{(\ell -1)k} a, 2^{(\ell -1)k+1} a, 2^{(\ell -1)k+2} a, \ldots, 2^{\ell k-1} a \right)
\]
is a geometric progression of length $k$ with common ratio 2.
If $A \in GPF_k(n)$, then $A$ does not contain this geometric progression, and so
at least one element in the set 
\[
X_{\ell}(a) = \left\{ 2^{(\ell -1)k} a, 2^{(\ell -1)k+1} a, 2^{(\ell -1)k+2} a, \ldots, 2^{\ell k-1} a \right\}
\]
 is not an element of $A$.
 Because every nonzero integer has a unique representation as the 
 product of an odd integer and a power of 2, it follows that, 
for  integers $\ell = 1,\ldots, L$ and  odd positive integers
$a \leq 2^{1-\ell k} n$, the sets $X_{\ell}(a)$
are pairwise disjoint subsets of $\{1,2,\ldots, n \}$.

For every real number $t \geq 1$, the number of odd positive integers 
not exceeding $t$ is strictly greater than $(t-1)/2$.
It follows that the cardinality of the set  $\{1,2,\ldots,n\} \setminus A$ 
is strictly greater than
\begin{align*}
\sum_{\ell = 1}^L \frac{1}{2} \left(  \frac{n}{2^{\ell k -1}} -1\right)
& =  \sum_{\ell = 1}^L \left(  \frac{n}{2^{\ell k}} - \frac{1}{2}\right) \\
& = n \sum_{\ell = 1}^L \frac{1}{2^{\ell k}} + O\left(\frac{\log n}{k} \right) \\
& = \frac{n}{2^k -1} + O\left(\frac{\log n}{k} \right).
\end{align*}

Note that if $r$ is an odd integer and $r \in X_{\ell}(a)$, then $\ell = 1$ and $r=a$.

Let $b$ be an odd  integer such that 
\beq  \label{GPF:ineq}
\frac{n}{6^{k-1}} < b \leq \frac{n}{5^{k-1}}
\eeq
and  $b$ is not divisible by 5, that is, 
\beq  \label{GPF:mod10}
b \equiv 1, 3,7,\text{ or } 9 \pmod{10}.
\eeq
We consider the following  geometric progression of length $k$ with ratio $5/3$:
\[
( 3^{k-1}b, 3^{k-2}5b,  \ldots,  3^{k-1-i}5^i b, \cdots, 5^{k-1} b ).
\]
Every integer in this progression is odd, and 
\[
\frac{n}{2^{k-1}} < 3^{k-1}b < \cdots < 5^{k-1}b \leq n.
\]
Let
\[
Y(b) = \{3^{k-1}b, 3^{k-2}5b,  \ldots,  3^{k-1-i}5^i b, \cdots, 5^{k-1} b \}.
\]
It follows that  
$X_{\ell}(a) \cap Y(b) = \emptyset$ for all $\ell$, $a$, and $b$.  
If the integers $b$ and $b'$ satisfy~\eqref{GPF:ineq} 
and~\eqref{GPF:mod10} with $b < b'$ 
and if $Y(b) \cap Y(b') \neq \emptyset$, then there exist
integers $i,j \in \{0,1,2,\ldots, k-1\}$ such that 
\[
3^{k-1-i}5^i b = 3^{k-1-j}5^j b'
\]
or, equivalently, 
\[
5^{i-j}b = 3^{i-j}b'.
\]
The inequality $b < b'$ implies that $0 \leq j < i \leq k-1$ and so
 $b' \equiv 0 \pmod{5}$, which contradicts~\eqref{GPF:mod10}.
Therefore, the sets $Y(b)$ are pairwise disjoint.
The number of integers $b$ satisfying inequality~\eqref{GPF:ineq} 
and congruence~\eqref{GPF:mod10}
is 
\[
 \frac{2}{5}\left( \frac{1}{5^{k-1} } -  \frac{1}{6^{k-1}} \right) n + O(1).
\]

Let $c$ be an odd integer such that 
\beq  \label{GPF:ineq-2}
\frac{n}{10^{k-1}} < c \leq \frac{n}{7^{k-1}}
\eeq
and $c$ is not divisible by 3 or 5, that is, 
\beq  \label{GPF:mod30}
c \equiv 1,7, 11, 13, 17, 19, 23, \text{ or } 29 \pmod{30}.
\eeq
We consider the following  geometric progression of length $k$ with ratio $7/5$:
\[
( 5^{k-1}c, 5^{k-2}7c,  \ldots,  5^{k-1-i}7^i c, \cdots, 7^{k-1} c ).
\]
Every integer in this progression is odd, and 
\[
\frac{n}{2^{k-1}} < 5^{k-1}c < \cdots < 7^{k-1}c \leq n.
\]
Let
\[
Z(c) = \{5^{k-1}c, 5^{k-2}7c,  \ldots,  5^{k-1-i}7^i c, \cdots, 7^{k-1} c \}.
\]
It follows that  
$X_{\ell}(a) \cap Z(c) = \emptyset$ for all $\ell$, $a$, and $c$.  
If $c$ and $c'$ satisfy~\eqref{GPF:ineq-2} and~\eqref{GPF:mod30} 
with $c < c'$  and if $Z(c) \cap Z(c') \neq \emptyset$, 
then there exist integers $i,j \in \{0,1,2,\ldots, k-1\}$ such that 
\[
5^{k-1-i}7^i c = 5^{k-1-j}7^j c'
\]
or, equivalently, 
\[
7^{i-j}c = 5^{i-j}c'.
\]
The inequality $c < c'$ implies that $0 \leq j < i \leq k-1$ and so
$c \equiv 0 \pmod{5}$, which contradicts~\eqref{GPF:mod30}.
Therefore, the sets $Z(c)$ are pairwise disjoint.

If $b$ and $c$ satisfy inequalities~\eqref{GPF:ineq} 
and~\eqref{GPF:ineq-2}, respectively, then $c < b$.  
If $Y(b) \cap Z(c) \neq \emptyset$, then there exist integers $i,j \in \{0,1,\ldots, k-1\}$
such that
\[
5^{k-1-i}7^i c = 5^{k-1-j}3^j b
\]
or, equivalently,
\[
5^j 7^i c = 5^i3^j b.
\]
Because $bc \not\equiv 0 \pmod{5}$, it follows that $i=j$
and so
\[
7^i c = 3^i b.
\]
Because $c < b$, we must have $i \geq 1$ and so $c \equiv 0 \pmod{3}$, 
which contradicts congruence~\eqref{GPF:mod30}. 
Therefore, $Y(b) \cap Z(c) = \emptyset$ and the sets $X_{\ell}(a)$, 
$Y(b)$, and $Z(c)$ are pairwise disjoint.  
The number of integers $c$ satisfying inequality~\eqref{GPF:ineq-2} 
and congruence~\eqref{GPF:mod30}
is 
\[
 \frac{4}{15}\left( \frac{1}{7^{k-1} } -  \frac{1}{10^{k-1}} \right) n + O(1).
\]
Because $A$ contains no $k$-term geometic progression, 
at least one element from each of the sets $X_{\ell}(a)$, 
$Y(b)$, and $Z(c)$ is not in $A$.
This completes the proof.
\end{proof}

\bc
If $A_k$ is a set of positive integers that contains no $k$-term geometric progression, then 
\[
d_U(A_k) \leq 1 - \frac{1}{2^k -1} - \frac{2}{5}\left( \frac{1}{5^{k-1} } -  \frac{1}{6^{k-1}} \right)
- \frac{4}{15}\left( \frac{1}{7^{k-1} } -  \frac{1}{10^{k-1}} \right).
\]
\ec

Here is a table of  upper bounds for $d_U(A)$ for various values of $k$:\\
\begin{center}
\begin{tabular}{|c|c|c|c|c|c|c|c|}\hline
k & 3& 4& 5&  6 &7 & 10 & 17 \\ \hline
$d_U(A_k) \leq$ & 0.84948   &   0.93147 &  0.96733  &   0.98404  & 0.99211& 0.99902 & 0.99999  \\
\hline
\end{tabular}
\end{center}

\section{Open problems}
For  every integer $k \geq 3$, let $GPF_k$ denote the set of sets
of positive integers that contain no $k$-term geometric progression.
It would interesting to determine precisely 
\[
\sup\{ d_U(A):  A \in GPF_k  \}
\]
and 
\[
\sup\{  d(A):  \text{$A$ has asymptotic density and $A \in GPF_k$}\}. 
\]
In the special case $k = 3$, Riddell~\cite[p. 145]{ridd69} claimed that if 
$A \in GPF_3$, then $d_U(A) < 0.8339$,
but wrote, "The details are too lengthy to be included here."

An infinite sequence  $A= (a_i)_{i=1}^{\infty}$ of positive integers  
is \emph{syndetic} if it is strictly increasing with bounded gaps.  
Equivalently, $A$ is syndetic if there is a number $c$ such that 
$1 \leq a_{i+1}-a_i \leq c$ for all positive integers $i$.  
Beiglb{\"o}ck, Bergelson, Hindman, and Strauss~\cite{beig-berg-hind-stra06} 
asked  if every syndetic sequence must contain arbitrarily long 
finite geometric progressions.

\def\cprime{$'$} \def\cprime{$'$} \def\cprime{$'$}
\providecommand{\bysame}{\leavevmode\hbox to3em{\hrulefill}\thinspace}
\providecommand{\MR}{\relax\ifhmode\unskip\space\fi MR }
\providecommand{\MRhref}[2]{%
  \href{http://www.ams.org/mathscinet-getitem?mr=#1}{#2}
}
\providecommand{\href}[2]{#2}

\end{document}